\documentclass[12pt]{article}
\usepackage{ytableau}
\usepackage{amssymb}
\usepackage[margin=1in]{geometry}
\usepackage{mathrsfs}
\usepackage{stmaryrd}
\usepackage{amsfonts,amsmath,amssymb,amscd}
\usepackage{shadow}
\usepackage{graphicx}
\usepackage{color}
\usepackage{longtable}
\usepackage{pstricks,multido}
\usepackage{hyperref}
\usepackage{qtree}
\allowdisplaybreaks

\newtheorem{thm}{Theorem}[section]
\newtheorem{lem}[thm]{Lemma}

\newtheorem{conj}{Conjecture}[section]

\def\qed{\hfill \rule{4pt}{7pt}}

\def\pf{\noindent {\it{Proof.} \hskip 2pt}}

\numberwithin{equation}{section}

\begin{document}
\begin{center}
{\large\bf The Raney Numbers and $(s,s+1)$-Core Partitions}
\end{center}

\begin{center}

{Robin D.P. Zhou}

College of Mathematics Physics and Information\\
ShaoXing University\\
ShaoXing 312000, P.R. China

zhoudapao@mail.nankai.edu.cn

\end{center}

\begin{abstract}
The Raney numbers $R_{p,r}(k)$ are a two-parameter generalization of the
Catalan numbers.
In this paper, we obtain a recurrence relation for the
Raney numbers which is a generalization of the recurrence relation
for the Catalan numbers.
Using this recurrence relation, we confirm a conjecture posed by Amdeberhan
concerning the enumeration of $(s,s+1)$-core partitions $\lambda$ with
parts that are multiples of $p$.
We then give a new combinatorial interpretation
for the Raney numbers $R_{p+1,r+1}(k)$ with $0\leq r<p$
in terms of $(kp+r,kp+r+1)$-core partitions $\lambda$ with
parts that are multiples of $p$.
\end{abstract}

\noindent {\bf Keywords}: Raney number, Catalan number, core partition, hook length, poset, order ideal, coral diagram

\noindent {\bf AMS  Subject Classifications}:  05A15, 05A17, 06A07

\section{Introduction}

In this paper, we build a connection between the Raney numbers and
 $(s,s+1)$-core partitions with
parts that are multiples of $p$.
We show that the number of $(kp+r,kp+r+1)$-core partitions
with parts that are multiples of $p$ equals the Raney number
$R_{p+1,r+1}(k)$, confirming a conjecture posed by Amdeberhan \cite{Amdeberhan}.

The Raney numbers were introduced by Raney in his investigation of
functional composition patterns \cite{Raney}
and these numbers  have also been used in
 probability theory \cite{Mlotkowski,Mlot}.
The Raney numbers are defined as follows:
\begin{equation}\label{exp:Raney}
R_{p,r}(k)=\frac{r}{kp+r}{kp+r \choose k}.
\end{equation}
The Raney numbers are a two-parameter generalization of the
Catalan numbers.
To be more specific, if $r=1$, the Raney numbers
specialize to the Fuss-Catalan numbers $C_p(k)$ \cite{Graham, Hilton},
where $C_p(k)$ are the numbers of  $k$-ary trees with labeled $p$
vertices and
\[C_p(k)=R_{p,1}(k)=\frac{1}{(p-1)k+1}{kp\choose k}.\]
If we further set $p=2$, we obtain the usual Catalan numbers $C_k$,
that is,
\begin{equation}\label{C=R}
R_{2,1}(k)=C_k=\frac{1}{k+1}{2k\choose k}.
\end{equation}

In \cite{Hilton}, Hilton and Pedersen obtained the following
relation for the Raney numbers $R_{p,r}(k)$ and the Fuss-Catalan numbers $C_p(k)$.

\begin{thm}\label{thm:Raney-Fuss}
Let p be a positive integer and let $r,k$ be nonnegative integers.
 Then we have
\begin{equation*}
R_{p,r}(k)=\sum_{i_1+\ldots+i_r=k}C_p(i_1)C_p(i_2)\ldots C_p(i_r).
\end{equation*}
\end{thm}

Based on Theorem \ref{thm:Raney-Fuss}, Beagley and Drube \cite{Beagley}
gave a combinatorial interpretation for the Raney numbers in terms of
coral diagrams which we will  describe in Section 2.

Let us give an overview of notation and terminology on partitions.
A partition $\lambda$ of a positive integer $n$ is a finite nonincreasing sequence of positive integers
$(\lambda_1, \lambda_2, \ldots, \lambda_m)$ such that
$\lambda_1+\lambda_2+\cdots+\lambda_m=n$.
We write
$\lambda=(\lambda_1,\lambda_2,\ldots,\lambda_m)\vdash n$
 and we say that $n$ is the size of $\lambda$ and $m$ is the
 length of $\lambda$.
The Young diagram of $\lambda$ is defined to be an up- and left-justified array
of $n$ boxes with $\lambda_i$ boxes in the $i$-th row.
Each box $B$ in $\lambda$ determines a hook
consisting of the box $B$ itself and
 boxes directly to the right and directly below $B$.
 The hook length of $B$, denoted $h(B)$, is the number
 of boxes in the hook of $B$.

For a partition $\lambda$,   the $\beta$-set of $\lambda$, denoted $\beta(\lambda)$, is defined to be the set of hook lengths of the boxes in the first column of $\lambda$.
For example, Figure
\ref{1.1fig} illustrates the Young diagram and the  hook lengths of a partition $\lambda=(5,3,2,2,1)$.
 The $\beta$-set of $\lambda$ is $\beta(\lambda)=\{9,6,4,3,1\}$.
Notice that a partition $\lambda$ is uniquely determined by its $\beta$-set.
Given a decreasing sequence of positive integers $(h_1,h_2,\ldots, h_m)$,
it is easily seen that the unique partition $\lambda$ with $\beta(\lambda)=\{h_1,h_2,\ldots, h_m\}$ is
\begin{equation}\label{betaset}
\lambda=(h_1-(m-1), h_2-(m-2), \dots, h_{m-1}-1, h_m).
\end{equation}

\begin{figure}[h]
\begin{center}
 \begin{ytableau}
     9& 7 & 4 & 2&1\\
     6 & 4 & 1\\
    4 & 2 \\
     3 & 1\\
     1
 \end{ytableau}
 \end{center}
\caption{The Young diagram of $\lambda=(5,3,2,2,1)$}
\label{1.1fig}
\end{figure}

For a positive integer $t$, a partition $\lambda$ is a $t$-core partition, or simply a $t$-core,  if
it contains no box whose hook length is a multiple of
 $t$.
Let $s$ be a positive integer not equal to $t$, we say that
$\lambda$ is an $(s,t)$-core if
it is simultaneously an $s$-core and a $t$-core.
For example, the partition $\lambda=(5,3,2,2,1)$ in
Figure \ref{1.1fig}  is a $(5,8)$-core.

Let $s$ and $t$ be two coprime positive integers.
Anderson \cite{Anderson} showed that the number of $(s,t)$-core partitions equals
${s+t\choose s}/(s+t)$.
It specializes to the Catalan number $C_s=\frac{1}{s+1}{2s \choose s}$
if $t=s+1$.
Ford, Mai and Sze \cite{Ford}  proved that the number of self-conjugate $(s,t)$-core partitions equals ${\lfloor \frac{s}{2}\rfloor+\lfloor\frac{t}{2}\rfloor \choose \lfloor\frac{s}{2}\rfloor}$.
Furthermore,   Olsson and  Stanton \cite{Olsson} proved
that there exists a unique $(s,t)$-core partition with the maximum size  $(s^2-1)(t^2-1)/24$.
 A simpler proof was
 provided by Tripathi \cite{Tripathi}.
More results on $(s,t)$-core partitions can be found
in  \cite{A-L,Armstrong,Chen,Johnson,Stanley,Wang,Zhou}.

In \cite{Amdeberhan},  Amdeberhan posed the following  conjecture.

\begin{conj} \label{conj}
 Let $s$ and $p$ be positive integers.
 The number of $(s,s+1)$-core partitions $\lambda$ with parts that are
 multiples of $p$ equals
\begin{equation}\label{equ:1.1}
\frac{s+1-p\lfloor\frac{s}{p}\rfloor}{s+1}{s+\lfloor\frac{s}{p}\rfloor\choose s}.
\end{equation}
\end{conj}

We observe that the expression (\ref{equ:1.1}) appearing in
Conjecture \ref{conj} equals the
 Raney number $R_{p+1,r+1}(k)$ if we write
$s=kp+r$, where $0\leq r < p$.
To prove conjecture \ref{conj}, in Section $2$, we shall deduce a
recurrence relation for the Raney numbers by using
coral diagrams.
In Section $3$, we shall give a characterization of the
$\beta$-set of the conjugate of an $(s,s+1)$-core partition
with parts that are multiples of $p$.
Based on this characterization,
we show that the number of $(kp+r,kp+r+1)$-core partitions
with parts that are multiples of $p$ have the same recurrence relation
with the Raney number $R_{p+1,r+1}(k)$.
This proves Conjecture \ref{conj}.

\section{Raney numbers}

In this section, we investigate the Raney numbers by using coral diagrams.
We obtain a recurrence relation for the Raney numbers which is a
generalization of the recurrence relation for the Catalan numbers.
We begin by introducing some graph theoretic terminology.

\begin{figure}[h]
\centerline{\includegraphics[width=15cm]{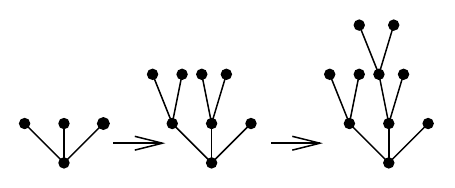}}
\caption{Construction of a coral diagram of type $(2,3,3)$}
\label{figcoral}
\end{figure}

Let $p$ be a positive integer. Then a p-star is a rooted tree with
$p$ terminal edges lying above a single base vertex.
A coral diagram of type $(p,r,k)$ is a rooted tree which is
constructed from an $r$-star via the repeated placement of
k $p$-stars atop terminal edges.
Let $D(p,r,k)$ denote the set of coral diagrams of type $(p,r,k)$.
We can construct a coral diagram
$D\in D(p,r,k)$ by attaching $p$-stars one ``tier'' at a time.
We begin with the base tree and work upward.
Figure \ref{figcoral} illustrates a coral diagram of type $(2,3,3)$.

We note that the definition of coral diagrams differs slightly from the
one defined by  Beagley and Drube \cite{Beagley}.
In \cite{Beagley}, a coral diagram of type $(p,r,k)$ is a rooted tree which is
constructed from an $(r+1)$-star via the repeated placement of
k $p$-stars atop terminal edges that are not the leftmost edge adjacent to
the root.
Not attaching $p$-stars to the leftmost edge adjacent to
the root gives them a consistent way of selecting
a base vertex for planar embedding.
Anyhow, the coral diagrams under the two definitions have the same enumeration.

We need the following theorem due to Beagley and Drube \cite{Beagley}.
\begin{thm}\label{thm:beagley}
Let $p$ be a positive integer and let $r,k$ be  nonnegative
integers.
Then the number of coral diagrams of type $(p,r,k)$ equals the
Raney number $R_{p,r}(k)$, that is,
\[|D(p,r,k)|=R_{p,r}(k).\]
\end{thm}

Using Theorem \ref{thm:beagley},
we obtain recurrence relations for Raney numbers, as follows.

\begin{thm}\label{thm:raney1}
Let $p$ be a positive integer and let $r,k$ be
 nonnegative integers.
Then we have
\begin{align}
 \label{equ:Raney1}
R_{p,1}(k)=&\sum_{i=0}^{k-1}R_{p,1}(i)R_{p,p-1}(k-1-i),\\[5pt]
\label{equ:Raney2}
 R_{p,r}(k)=&\sum_{i=0}^{k}R_{p,1}(i)R_{p,r-1}(k-i),  r>1.
\end{align}
And the Raney numbers $R_{p,r}(k)$ are uniquely determined by the above relations
 together with initial values $R_{p,r}(0)=1$, $R_{p,0}(k+1)=0$.
\end{thm}

\pf
The initial values $R_{p,r}(0)=1$, $R_{p,0}(k+1)=0$
follow directly from the
expression (\ref{exp:Raney}) for the Raney numbers, that is,
 \begin{equation*}
R_{p,r}(k)=\frac{r}{kp+r}{kp+r \choose k}.
\end{equation*}

Now we assume that $k>0$, $r>0$.
From Theorem \ref{thm:beagley},
we see that $R_{p,r}(k)$ equals the number of coral diagrams of type $(p,r,k)$,
that is, \[R_{p,r}(k)=|D(p,r,k)|.\]
We proceed to investigate the recurrence relations for Raney numbers $R_{p,r}(k)$
by using coral diagrams.
Let $D$ be a coral diagram in $D(p,r,k)$.
There are two cases.

\begin{figure}[h]
\centerline{\includegraphics[width=10cm]{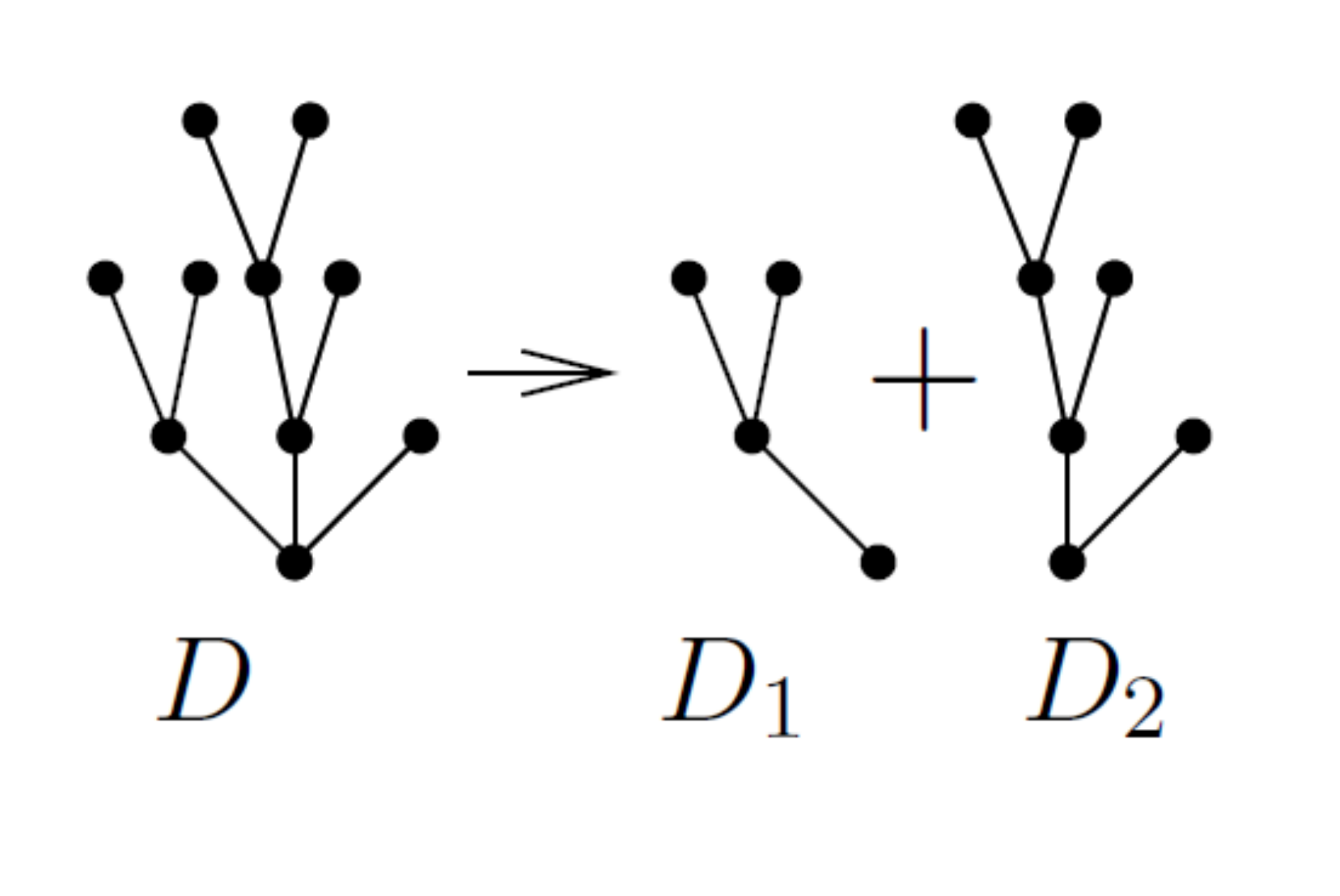}}
\caption{Dividing a coral diagram into two coral diagrams}
\label{figdividing}
\end{figure}

If $r>1$, we divide the coral diagram $D$ into
two coral diagrams $D_1$ and $D_2$ by splitting the root of $D$ such
that the root of the coral diagram $D_1$ is only adjacent
 to the leftmost edge and the root of the coral diagram $D_2$ is
adjacent to the rest $r-1$ edges.
See Figure \ref{figdividing} as an example.
Assume that the coral diagram $D_1$ has $i$ $p$-stars, then
 the coral diagram $D_2$ has $k-i$ $p$-stars.
It yields that $D_1 \in D(p,1,i)$ which is counted by $R_{p,1}(i)$,
and $D_2 \in D(p,r-1,k-i)$ which is counted by $R_{p,r-1}(k-i)$.
Since the coral diagram $D$ is uniquely determined by $D_1$ and $D_2$,
we obtain the relation (\ref{equ:Raney2}), that is,
\[R_{p,r}(k)=\sum_{i=0}^{k}R_{p,1}(i)R_{p,r-1}(k-i).\]

If $r=1$, then $D$ is constructed from a 1-star.
Contracting  the edge of this 1-star into a new vertex,
we obtain a new coral diagram $D'$.
It is easily seen that $D'$ is a coral diagram in
$D(p,p,k-1)$, which is counted by $R_{p,p}(k-1)$.
Hence we have $R_{p,1}(k)=R_{p,p}(k-1)$.
If further $p=1$, since $R_{1,0}(k)=0$ when $k>0$,
relation (\ref{equ:Raney1}) is equivalent to
$R_{1,1}(k)=R_{1,1}(k-1)$.
By the expression (\ref{exp:Raney}) for the Raney numbers $R_{p,r}(k)$,
it can be easily checked that $R_{1,1}(k)=R_{1,1}(k-1)=1$.
Hence relation (\ref{equ:Raney1}) holds if $r=p=1$.
Now assume that $p>1$.
Combining relation (\ref{equ:Raney2}) and $R_{p,1}(k)=R_{p,p}(k-1)$,
we obtain (\ref{equ:Raney1}).

To prove that the Raney numbers are uniquely determined by the relations
(\ref{equ:Raney1}) and (\ref{equ:Raney2})
 together with initial values $R_{p,r}(0)=1$, $R_{p,0}(k+1)=0$,
we use induction on $k$.
If $k=0$, we have that $R_{p,r}(0)=1$.
Suppose we can determine all the values of the Raney
numbers $R_{p,r}(k)$ with $k \leq n$ from the
relations
(\ref{equ:Raney1}) and (\ref{equ:Raney2})
 together with initial values.
We proceed to show that $R_{p,r}(n+1)$ can be determined.
First, $R_{p,1}(n+1)$ can be determined by  (\ref{equ:Raney1}).
Then by (\ref{equ:Raney2}), the values of
$ R_{p,2}(n+1), R_{p,3}(n+1),R_{p,4}(n+1)\cdots$ can be
obtained in turn.
Thus we can determine all the values of the Raney numbers
$R_{p,k}(n+1)$.
By induction, all the  Raney numbers $R_{p,r}(k)$ can be
determined.
This completes the proof.
\qed

From (\ref{C=R}), we have that $C_k =R_{2,1}(k)$.
Substituting $p=2$ into (\ref{equ:Raney1}),
we obtain the recurrence relation for the Catalan numbers $C_k=\sum_{i=0}^{k-1}C_iC_{k-1-i}$.

\section{Proof of Conjecture \ref{conj}}

In this section, we shall give a characterization of the
$\beta$-set of the conjugate of an $(s,s+1)$-core partition
with parts that are multiples of $p$.
Using this characterization, we prove Conjecture \ref{conj}
posed by Amdeberhan.
We then give a new combinatorial interpretation
for the Raney numbers $R_{p+1,r+1}(k)$ with $0\leq r<p$
in terms of $(kp+r,kp+r+1)$-core partitions $\lambda$ with
parts that are multiples of $p$.

We observe that the expression (\ref{equ:1.1}) appearing in
Conjecture \ref{conj}  equals the Raney number $R_{p+1,r+1}(k)$ if we write
$s=kp+r$, where $0\leq r < p$.
That is,
\begin{align*}
&\frac{s+1-p\lfloor\frac{s}{p}\rfloor}{s+1}{s+\lfloor\frac{s}{p}\rfloor\choose s}=
\frac{kp+r+1-kp}{kp+r+1}{kp+r+k\choose kp+r}  \\[6pt]
&=\frac{r+1}{k(p+1)+r+1}{k(p+1)+r+1\choose k}=R_{p+1,r+1}(k).
\end{align*}

Hence Conjecture \ref{conj} can be restated as the following equivalent theorem.

\begin{thm}\label{thm:main}
 Let $s$ and $p$ be positive integers.
Suppose that $s=kp+r$, where $0\leq r<p$.
Then the number of $(s,s+1)$-core partitions $\lambda$ with parts that are
 multiples of $p$ equals the Raney number $R_{p+1,r+1}(k)$.
\end{thm}

To prove the above theorem, we shall
give a characterization of the
$\beta$-set $\beta(\lambda^c)$,
where $\lambda^c$ is the conjugate of an $(s,s+1)$-core partition $\lambda$ with parts that are multiples of $p$.
Let us  recall some notation and terminology on posets.

Let $P$ be a poset. For two elements $x$ and $y$ in $P$, we say $y$ covers $x$ if $x< y$ and there exists no element $z\in P$
satisfying $x< z< y$.
 The Hasse diagram of a finite poset $P$ is a graph whose vertices are the elements of $P$, whose edges are the cover relations, and such that if $y$ covers $x$ then there is an edge connecting $x$ and $y$ and
 $y$ is placed above $x$.
 An order ideal of $P$ is
a subset $I$ such that  if any $y\in I$ and $x\leq y$ in $P$, then $x\in I$. Let $J(P)$ denote the set of order ideals of  $P$.
 For more details on poset, see Stanley \cite{Stanleybook}.

In the following lemma, Anderson \cite{Anderson} established a correspondence between core partitions and
 order ideals of a certain poset by mapping a partition to its $\beta$-set.

\begin{lem}\label{lem:core-poset}
Let $s,t$ be two coprime positive integers, and let $\lambda$ be a partition of $n$. Then $\lambda$ is an  $(s,t)$-core partition if and only if $\beta(\lambda)$ is an order ideal of $P_{(s,t)}$.
where
\[P_{(s,t)}=\mathbb{N}^+\setminus \{n \in\mathbb{N}^+\mid n=k_1s+k_2t \mbox{ for some } k_1,k_2\in \mathbb{N}\}\]
and $y$ covers $x$ in $P_{(s,t)}$ if  $y-x\in\{s,t\}$.
\end{lem}

\begin{figure}[h]
\centerline{\includegraphics[width=10cm]{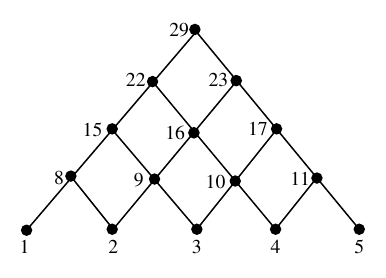}}
\caption{The Hasse diagram of $T_6$}
\label{Hasse}
\end{figure}

For example, let $s=3$ and $t=4$.
We can construct all $(3,4)$-core partitions by finding order ideals of $P_{(3,4)}$.
It is easily checked that $P_{(3,4)}=\{1,2,5\}$ with the
partial order $5>2$ and $5>1$. Hence the order ideals of $P_{(3,4)}$ are
$\emptyset$, $\{1\}$, $\{2\}$, $\{2,1\}$ and $\{5,2,1\}$.
The corresponding $(3,4)$-core partitions
are $\emptyset$, $(1)$, $(2)$, $(1,1)$ and $(3,1,1)$, respectively.
Let $T_s=P_{(s,s+1)}$.
From Lemma \ref{lem:core-poset}, the Hasse diagram of $T_s$ can be easily constructed.
For example, Figure \ref{Hasse} illustrates the Hasse diagram of the poset $T_6$.

Let $\lambda$ be an $(s,s+1)$-core partition with parts that are
 multiples of $p$ and let
$\lambda^c$ be the conjugate of $\lambda$.
Since $\lambda^c$ is also an $(s,s+1)$-core partition,
from Lemma \ref{lem:core-poset}, we see that
$\beta(\lambda^c)$ is an order ideal in $T_s$.
The following lemma gives a further
description of $\beta(\lambda^c)$.

\begin{lem}\label{lem}
Let $\lambda$ be a partition  with parts that are
multiples of $p$ and let $\lambda^c$ be the conjugate of $\lambda$.
Then the length of each maximal consecutive integers in $\beta(\lambda^c)$
is divided by $p$.
\end{lem}

\pf
Suppose that
$\lambda^c=(a_1^{m_1},a_2^{m_2},\ldots a_n^{m_n})$, where
$a_1>a_2> \cdots>a_n$ and $a_i^{m_i}$ means $m_i$ occurrences of $a_i$.
Since $\lambda$ is a partition  with parts that are
multiples of $p$,
we have that $m_i(1\leq i \leq n)$ is a multiple of $p$.
Then the lemma follows directly from the relation of a partition and
its $\beta$-set shown in (\ref{betaset}).
\qed

Let $S$ be an integer set and let $l$ be a positive integer.
We say that the set  $S$ has property $P_l$
if the length of each maximal consecutive integers in $S$
is divided by $l$.
For example, the integer set $\{1,2,3,6,7,8,9,10,11\}$ has
property $P_3$.

Let $\lambda$ be a partition.
From Lemma \ref{lem:core-poset} and Lemma \ref{lem},
we see that $\lambda$ is an $(s,s+1)$-core partition with parts that are
multiple of $p$ if and only if $\beta(\lambda^c)$ is an order ideal in $T_s$
with property $P_p$.
Let $s$ and $p$ be positive integers.
Suppose that $s=kp+r$, where $0\leq r<p$.
Let $C_{p,r}(k)$ denote  the number of order ideals in $T_{s}$ with
property $P_p$.
We set $C_{p,-1}(k)=0$ if $k>0$ and $p\geq 0$.
Hence to prove Theorem \ref{thm:main},
 it suffices to show that  $C_{p,r}(k)=R_{p+1,r+1}(k)$.

We proceed to compute the number of  order ideals in $T_{s}$ with
property $P_p$.
To this end, we shall partition $J(T_s)$ according to the smallest missing element of rank $0$ in an order ideal.
Note that the elements of rank $0$ in $T_s$ are just the minimal elements.
For $1\leq i\leq s-1$, let $J_i(T_s)$ denote the set of order ideals of $T_s$ such that $i$ is the smallest missing element of rank $0$.
Let $J_s(T_s)$ denote the set of order ideals which contain all minimal elements in $T_s$.
Then we can write $J(T_s)$ as \[ J(T_s)=\bigcup_{i=1}^sJ_i(T_s).\]

\begin{figure}[h]
\centerline{\includegraphics[width=15cm]{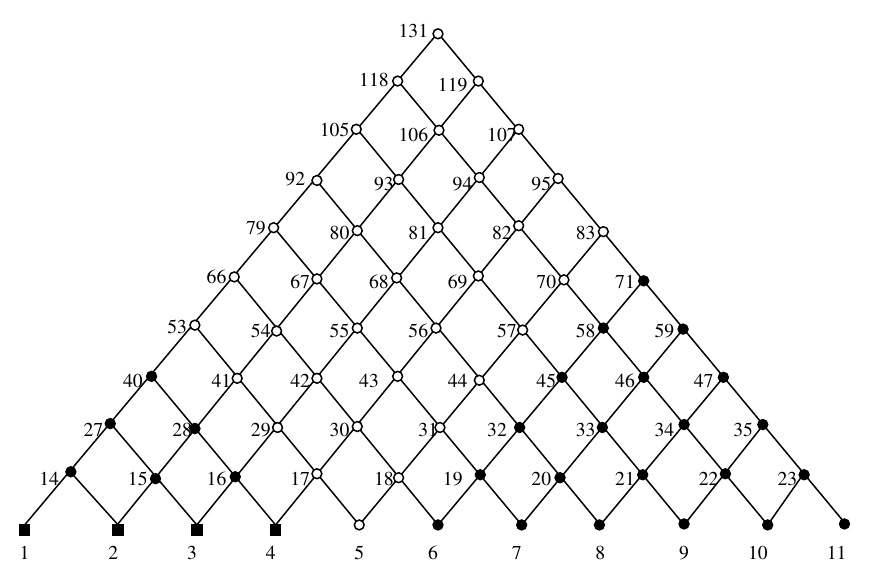}}
\caption{The  elements of an order ideal $I\in J_5(T_{12})$}
\label{Hasse1}
\end{figure}

Figure \ref{Hasse1} gives an illustration of the   elements contained in an order ideal in $J_5(T_{12})$.
We see that an order ideal $I\in J_5(T_{12})$ must contain the elements labeled by squares, but does not contain any elements represented by open  circles.
The elements represented by solid  circles  may or may not  appear in $I$.
That is, $I$ can be decomposed into three parts, one is $\{1,2,3,4\}$, one is isomorphic to an order ideal of $T_4$ and one is isomorphic to an order ideal of $T_7$.

In general, for $1\leq i\leq s$ and an order ideal $I\in J_i(T_s)$, we can decompose it into three parts: one is $\{1,2,\ldots,i-1\}$, one is isomorphic to an order ideal of $T_{i-1}$ and one is isomorphic to an order ideal of $T_{s-i}$
(Some parts may be empty).
We shall use  this decomposition to prove Theorem \ref{thm:main}.

\noindent{\it Proof of Theorem \ref{thm:main}.}
We prove this theorem by
 showing that $C_{p,r}(k)$ and $R_{p+1,r+1}(k)$ have
the same recurrence relation and initial values.
If $k=0$, that is, $s=r<p$.
Then the unique order ideal in $T_{r}$ with
property $P_p$ is $\emptyset$.
By the definition of $C_{p,r}(k)$,
we have that $C_{p,r}(0)=1=R_{p+1,r+1}(0)$.
Combining with $C_{p-1,-1}(k)=R_{p,0}(k)=0$, where $k>0$,
we see that $C_{p,r}(k)$ have the same initial values with $R_{p+1,r+1}(k)$.

We proceed to show that $C_{p,r}(k)$ and $R_{p+1,r+1}(k)$ have
the same recurrence relation.
Suppose now $k\geq 1$.
Let $I$ be an order ideal in $J(T_{s})$ with property $P_p$.
Then $I\in J_{ip+1}(T_s)$, where $0\leq i\leq  \lfloor\frac{s-1}{p}\rfloor $ since
$I$ has property $P_p$.
$I$ can be decomposed into three parts:
one is $\{1,2,\ldots,ip\}$, one is isomorphic to an order ideal $I_1$ of $T_{ip}$ and one is isomorphic to an order ideal $I_2$ of $T_{s-1-ip}$.
Since the absolute difference of any two numbers in two parts are larger than $1$,
we have that all of the three parts $\{1,2,\ldots,ip\}$, $I_1$ and $I_2$ have
property $P_p$.

It is easily seen that the number of order ideals in $T_{ip}$
with property $P_p$ is counted by $C_{p,0}(i)$.
To enumerate the number of order ideals in $T_{s-1-ip}$
with property $P_p$,
we consider two cases.
If $r=0$, namely $s=kp$, then $s-1-ip=(k-i-1)p+(p-1)$.
It follows that the number of order ideals in $T_{s-1-ip}$
with property $P_p$ is counted by $C_{p,p-1}(k-i-1)$.
If $r>0$, then $s-1-ip=(k-i)p+r-1$.
So in this case, the number of order ideals in $T_{s-1-ip}$
with property $P_p$ is counted by $C_{p,r-1}(k-i)$.

Combining the two cases, we have that
\begin{align*}
C_{p,0}(k)=&\sum_{i=0}^{k-1}C_{p,0}(i)C_{p,p-1}(k-i-1),\\[5pt]
 C_{p,r}(k)=&\sum_{i=0}^{k}C_{p,0}(i)C_{p,r-1}(k-i),  r>0.
\end{align*}
It can be checked that $C_{p,r}(k)$  have the same recurrence
relation with the Raney numbers
$R_{p+1,r+1}(k)$ as shown in Theorem \ref{thm:raney1}.
By a similar discussion as in the proof of
Theorem \ref{thm:raney1},
we obtain that $C_{p,r}(k)$ are uniquely determined by the
above relations together with the initial values.
Hence we have that $C_{p,r}(k)=R_{p+1,r+1}(k)$.
This completes the proof.
\qed

From Theorem \ref{thm:main}, we see that the Raney numbers
$R_{p+1,r+1}(k)$ with $0\leq r<p$  equal the numbers of
$(kp+r,kp+r+1)$-core partitions with parts that are
multiples of $p$.
This gives a new combinatorial interpretation for
these Raney numbers.

\vspace{0.5cm}
 \noindent{\bf Acknowledgments.}
This work was supported by the National Science Foundation of China.

\end{document}